# Stochastic Real-time Power Dispatch with Large-scale Wind Power Integration and its Analytical solution

Shuwei Xu, Wenchuan Wu, Yue Yang, Bin Wang, Xiaohai Wang

*Abstract*—Real-time power dispatch (RTD) can coordinate wind farms, automatic generation control (AGC) units and non-AGC units. In RTD, the probable wind power forecast errors (WPFE) should be appropriately formulated to ensure system security with high probability and minimize operational cost. Previous studies and our onsite tests show that Cauchy distribution (CD) effectively fits the "leptokurtic" feature of small-timescale WPFE distributions. In this paper, we propose a chance-constrained real-time dispatch (CCRTD) model with the WPFE represented by CD. Since the CD is stable and has promising mathematical characteristics, the proposed CCRTD model can be analytically transformed to a convex optimization problem considering the dependence among wind farms' outputs. Moreover, the proposed model incorporates an affine control strategy compatible with AGC systems. This strategy makes the CCRTD adaptively take into account both the additional power ramping requirement and power variation on transmission lines caused by WPFE in RTD stage. Numerical test results show that the proposed method is reliable and effective. Meanwhile it is very efficient and suitable for real-time application.

*Index Terms*-- Stochastic optimization, real-time dispatch, wind power forecast error, Cauchy distribution

## Nomenclature

### A. Indices

| | |
|---|---|
| $t$ | Indices for time periods. |
| $i$ | Indices for conventional units (non-AGC units). |
| $j$ | Indices for AGC units. |
| $k$ | Indices for wind farms. |
| $d$ | Indices for loads. |
| $l$ | Indices for transmission lines. |

### B. Parameters and Functions

| | |
|---|---|
| $\mathbf{x}$ | Random vector |
| $x$ | The value of random vector/variable |
| $\boldsymbol{\mu}, \boldsymbol{\Sigma}$ | Location vector and scale matrix of multivariate CD |
| $\mu, \sigma^2$ | Location and scale parameters of one-dimensional CD |
| $T$ | Number of time periods. |
| $\Delta T$ | Length of each time period in minutes. |
| $N$ | Number of conventional units. |
| $J$ | Number of AGC units. |
| $K$ | Number of wind farms. |
| $D$ | Number of loads / nodes. |
| $L$ | Number of transmission lines. |
| $\bar{p}_{k,t}^w$ | The maximal power output for wind farm $k$ during period $t$ (MW). |
| $\bar{w}_{j,t}$ | Upper bound of total available wind power during period $t$ (MW). |
| $\bar{P}_{i,t}^s / \underline{P}_{i,t}^s$ | Upper / lower limitation of active power output of non-AGC unit $i$ during period $t$ (MW). |
| $\bar{P}_{j,t}^a / \underline{P}_{j,t}^a$ | Upper / lower limitation of active power output of AGC unit $j$ during period $t$ (MW). |
| $\alpha_{j,t}$ | Participation factor of AGC unit $j$ during time period $t$. |
| $CF_{i,t}(\cdot)$ | Generation cost of conventional unit $i$ during period $t$ ($). |
| $CF_{j,t}(\cdot)$ | Generation cost of AGC unit $j$ during period $t$ ($) |
| $CR_{j,t}^+(\cdot) / CR_{j,t}^-(\cdot)$ | Upward / downward regulation cost of AGC unit $j$ during period $t$ ($). |
| $E(\cdot)$ | Expectation of random variables. |
| $a_{i,t} / b_{i,t} / c_{i,t}$ | Cost coefficient of conventional unit $i$ in period $t$. |
| $a_{j,t} / b_{j,t} / c_{j,t}$ | Cost coefficient of AGC unit $j$ in period $t$. |
| $\varphi_t(\cdot)$ | Probability density function of random variable during time period $t$. |
| $\Phi_t(\cdot)$ | Cumulative distribution function of random variable during time period $t$. |
| $\gamma_{j,t}^+ / \gamma_{j,t}^-$ | Coefficient of the upward / downward regulation cost of AGC unit $j$ during period $t$ ($/MW). |
| $p_{d,t}^d$ | the load demands of node $d$ in period $t$. (MW). |
| $RU_{i,t}^s / RD_{i,t}^s$ | Upward / downward ramping rate of non-AGC unit $i$ during period $t$ (MW/min). |
| $RU_{j,t}^a / RD_{j,t}^a$ | Upward / downward ramping rate of AGC unit $j$ during period $t$ (MW/min). |
| $R_t^+ / R_t^-$ | Upward / downward reserve requirement during period $t$ (MW). |
| $G_{l,i} / G_{l,k} /$ $G_{l,j} / G_{l,d}$ | Power transfer distribution factors of line $l$ with respect to $p_{i,t}^s, \tilde{p}_{j,t}^a, \tilde{p}_{k,t}^w, p_{d,t}^d$. |
| $L_{l,t}$ | Power flow limits for transmission line $l$ during period $t$ (MW). |

Manuscript received XX, 2019. This work was supported in part by the National Key R&D Program of China (2018YFB0904200).

S. Xu, W. Wu, Y. Yang and B. Wang are with the State Key Laboratory of Power Systems, Department of Electrical Engineering, Tsinghua University, Beijing 100084, China (e-mail: wuwench@tsinghua.edu.cn). X. Wang is Inner Mongolia Power (Group) Co, Ltd.

| $\delta, \beta, \varepsilon, \eta$ | pre-defined allowed probabilities of violation |

## C. Deterministic Variables

| $p_{k,t}^w$ | Scheduled power output of wind farm $k$ during period $t$ (MW). |
| $w_t$ | The summation of scheduled total wind power during period $t$ (MW). |
| $p_{j,t}^a$ | Scheduled base-point of AGC unit $j$ during period $t$ (MW). |
| $p_{i,t}^s$ | Scheduled active power output of non-AGC unit $i$ during period $t$ (MW). |

## D. Uncertain Variables

| $\tilde{p}_{k,t}^w$ | Actual wind power output of wind farm $k$ during period $t$ (MW). |
| $\tilde{w}_t$ | The summation of actual total wind power output during period $t$ (MW). |
| $\tilde{p}_{j,t}^a$ | Actual power output of AGC unit $j$ during period $t$ (MW). |

## I. INTRODUCTION

### A. Background and challenge

WIND power penetration has increased significantly in recent years [1]. Large-scale wind power integration can provide green energy and bring environmental benefits [2]. However, because of the intermittency and randomness of wind power [2], system operators need to consider the inherent challenges in real-time dispatch(RTD). RTD is essentially formulated as a dynamic economic dispatch model (ED) in which only the solution of the first period is committed. Stochastic ED models are conventionally used to deal with the uncertainties of wind power but with two crucial issues. The first is how to describe the forecast error precisely and model-friendly. The second is how to incorporate uncertainties into operational constraints and objective function properly to formulate a tractable stochastic ED model. A detailed description of these two challenges is as follows.

a) To reduce the operation costs and increase the system reliability, an "accurate" WPFE model with smaller fitting error is necessary [3] first. And then, due to similar or discrepant meteorological conditions, wind power of multiple wind farms exhibits correlation or complementary characteristics on different regions [4] which should be contained in an " accurate" WPFE model. Furthermore, to facilitate the operation of power system with large-scale wind power penetration, a "model-friendly" WPFE model usually has these mathematical properties: (1) The linear combination of related random variables can be expressed by same type distribution (this is the concept of "stable distribution"[5]) or a simple distribution; (2)The inverse cumulative distribution function(*CDF*) should be expressed analytically in order to convert chance constraints into deterministic constraints. (3) The expected cost in the form of integral can be given analytically. WPFE models with part of the above properties have shown advantages in ED problems [6][7]. But none of them have ability to access all these properties simultaneously.

b) For the process of uncertainty modeling and problem solving, chance-constrained economic dispatch (CCED) [8] with adjustable confidence level [9] is a good choice for balancing the security and economy of the dispatch process. However, because of the fluctuation of wind power and the redistribution of unbalance power between regulation generators, the branch flow and power ramping are all stochastic in the operational constraints, which makes CCED more complicated. Moreover, the tightest bottleneck in solving chance-constrained optimization problems is inefficiency, which prevents its application in real-time dispatch.

### B. Literature review

*1) Wind power forecast model*

WPFEs are traditionally assumed as random variables [1] that follow certain distribution. Some studies have shown that the feature of small-timescale WPFE distributions is "leptokurtic" [10], which means the property of both high kurtosis [1][11] and fat-tail [12]. "High kurtosis" indicates high prediction accuracy while "fat tails" provides the frequency of extreme events [10] and holds important information in reliability studies [13]. In paper [14], the author examined the errors from operational wind power forecasting systems, finding that WPFE were poorly represented by normal distribution [15]. Some other existing models, such as Weibull [16] or beta distribution [17] were also not suited to fully describe the heavy-tailed character of WPFE data [13]. For power systems with multiple wind farms, comprehensive analyzing the output correlation and dependence between various wind farms is necessary in dynamic ED [18][19]. Neglecting the correlation can lead extra cost and increase the risk of transmission lines overloading [4]. In literature [20], Xie proposed a novel data-driven wind speed forecast framework by leveraging the spatial-temporal correlation among geographically dispersed wind farms and then wind power forecasts were converted from wind speed forecasts based on power curve. Other studies used copula function[18] [21] [22] to formulate the correlation or dependence of multi-wind farms generation. But scenario sampling brings heavy computational burden, meanwhile, it is hard to select a suitable or optimal copula function. Besides, the multivariate joint distribution function, such as multivariate Gaussian distribution [8], was also used to formulate the correlation.

As a kind of stable distribution, Cauchy distribution (CD) is good at modeling the high spike of the frequency histogram [11] with fat tails [10]. In previous study [1], the forecast error distribution was fitted using maximum-likelihood optimization, and it was found that the CD was better than the Gaussian, Beta and Weibull distributions in all cases. Moreover, the dependence of all wind farms' output can be described by the scale matrix of multivariate CD.

*2) Stochastic ED model and solution*

To capture the uncertainty and fluctuation of wind power, stochastic optimization (SO) models are widely used for ED. Stochastic optimization model with chance constraints allows to trade off security for economy with adjustable risk level to meet the different reliability requirements [8]. But it is hard to solve even if the problem is convex [23]. In [23], a scenario-based approach was developed to replace chance constraints by sampling the uncertainty parameters which were pre-processed in an offline stage to accelerate the real-time decision in online stage. Reference [24] proposed bilinear and linear formulations for the chance-constrained mixed integer programming with Benders decomposition method based on Monte Carlo simulation. In paper [25], the authors improved the stochastic programming approach in the unit commitment problem to incorporate wind power scenarios by introducing a dynamic decision making approach. The expectation in objective function was approximately calculated by scenario sampling. In paper [26][27], a chance-constrained two stage stochastic program was solved by a combined sample average approximation(SAA) algorithm. It is found that all of these methods are based on scenario generation. If the sample size is large enough, feasibility in the chance-constrained sense can be guaranteed with high confidence [23][25]. However, scenario-based method suffers heavy computational burden which limits its application to real-time decision making.

The versatile distribution(VD) and Truncated Versatile Distribution(TVD) in paper [6] and [28] were proposed to model the WPFE. Compared with the Gaussian and Beta distributions, the VD/TVD can fit the WPFE more accurately and their inverse *CDFs* have analytical mathematical expressions. VD/TVD-based chance constraints can be converted into deterministic constraints with quantile. However, it cannot provide a distribution for a linear combination of random variables represented by VD/TVD, so the transmission capacity constraints were ignored in there works. Recently, the Gaussian mixture model (GMM) was employed to model the correlated prediction errors of wind power in CCED [7] or ramping capacity allocation [9]. In these works, the *CDF* of a Gaussian distribution was proximately fitted by a piecewise fourth-order polynomial and then the chance constraints were converted into deterministic constraints. The drawbacks of these formulations include: (1) the piecewise characteristic of each polynomial component mitigates against a direct solution of the inverse *CDF*; and, (2) the existence of multiple solutions to the quartic equation means that an extra validation step is required to obtain a reasonable solution.

*C. Contributions*

In this paper, we statistically analyze the small-timescale WPFE distributions for 20 wind farms in Southwest China and the onsite tests justify that the CD outperforms the other distributions, such as Gaussian, Beta or Weibull distributions, especially in capturing the kurtosis and tail behavior, which also has been demonstrate in paper [1] and [10]. Therefore, we apply the CD to characterize uncertainties of WPFE in the CCRTD model with affine AGC control strategy (A-CCRTD). Due to the promising characteristics of CD, the A-CCRTD is analytically converted into a deterministic convex optimization problem and solved efficiently without any approximation. In more detail, the contributions of this paper include:

1) Since the remarkable feature of WPFE is "leptokurtic" and the effectiveness of modeling the observed WPFE with a CD over short timescales, we formulate the wind farms' output as a multivariate random variable represented by a multivariate CD. The dependence between multiple wind farms' output is also described in the scale matrix. We are the first to apply the CD in characterizing the WPFE in CCED.

2) As a kind of analytical "stable distribution", CD has several favorable mathematical properties: its *CDF* and inverse *CDF* can be expressed analytically; the expected cost in the form of integral with CD can also be given analytically. Therefore, we can easily convert all the chance linear constraints into deterministic linear constraints and the expectation terms in the objective function into analytical convex functions precisely. Finally, the A-CCRTD can be solved efficiently with a guaranteed global optimal solution. This approach makes this chance-constrained optimization problem practically tractable in real-time applications for power systems with high wind power penetration. Besides, since the transformation of chance constraints is analytical, sensitivity analysis such as tuning risk level can be achieved easily.

3) An affine control strategy for the AGC system is incorporated into the chance-constrained dispatch process. The proposed model takes into account both the additional power ramping requirement (APRR) and the power variation on transmission lines caused by the WPFEs in RTD stage, i.e., sufficient regulation capacity of AGC units should be reserved in RTD stage to correct the real-time power mismatch caused by wind power uncertainties to ensure system reliability.

The remainder of this paper is organized as follows. In Section II, the mathematical properties of Cauchy distribution and WPFE model are presented. Section III is the mathematical formulation of the A-CCRTD model with Cauchy distribution. Section IV provides case studies and Section V is conclusion.

## II. WPFE Modeling With Multivariate Cauchy Distribution

If a *p*-dimensional random vector **x** follows a multivariate Cauchy distribution with location vector $\boldsymbol{\mu}$ and scale matrix $\Sigma$, that is $\mathbf{X} \sim Cauchy_p(\boldsymbol{\mu}, \Sigma)$, then the *PDF* is presented as [29]:

$$f_X(x;\boldsymbol{\mu},\Sigma) = \frac{\Gamma\left(\frac{1+p}{2}\right)}{\Gamma\left(\frac{1}{2}\right)\pi^{\frac{p}{2}}|\Sigma|^{\frac{1}{2}}\left[1+(x-\boldsymbol{\mu})^T\Sigma^{-1}(x-\boldsymbol{\mu})\right]^{\frac{1+p}{2}}} \quad (1)$$

where $x, \boldsymbol{\mu} \in \mathbb{R}^p$ and $\Sigma \in \mathbb{R}^{p \times p}$ is a positive-definite matrix.

When p = 1, that is $\mathbf{X} \sim Cauchy(\mu, \sigma^2)$, the probability density function (*PDF*) of one-dimensional Cauchy distribution is:

$$f_x(x;\mu,\sigma^2) = \frac{1}{\pi}\left[\frac{\sigma}{(x-\mu)^2+\sigma^2}\right], \quad x \in \mathbb{R} \quad (2)$$

Some important properties beneficial to solve A-CCRTD are listed below:

*1) Integral property:*

$$\int x \cdot PDF(x)dx = \frac{\sigma}{2\pi}\ln(1+(\frac{x-\mu}{\sigma})^2) + \frac{\mu}{\pi}\arctan\left(\frac{x-\mu}{\sigma}\right) + c \quad (3)$$

$$\int x^2 \cdot PDF(x)dx = \frac{\sigma}{\pi}(x-\mu) + \frac{(\mu^2-\sigma^2)}{\pi}\arctan\left(\frac{x-\mu}{\sigma}\right) \\ + \frac{\mu\sigma}{\pi}\ln(1+(\frac{x-\mu}{\sigma})^2) + c \quad (4)$$

*2) Stable property:*

"Stable" [5] means the linear transformation of *x* in equation (1) can be expressed as a new random variable comply with one-dimensional Cauchy distribution. For example, suppose that *a* is a *p*-dimensional vector, and then we have

$$a^T x \sim Cauchy\left(a^T \boldsymbol{\mu}, a^T \Sigma\, a\right) \quad (5)$$

## 3) Analytical expressions of CDF and inverse CDF

$$CDF(x) = \frac{1}{\pi}\arctan\left(\frac{x-\mu}{\sigma}\right) + \frac{1}{2} \quad (6)$$

$$CDF^{-1}(F) = \mu + \sigma \tan\left[\pi\left(F - \frac{1}{2}\right)\right] \quad (7)$$

where $F$ is the quantile.

## 4) Fitting and sampling

The multivariate CD parameters are fitted to the data using the mscFit function of the fMultivar package [30] in the R statistical computing environment [31]. While, the sampling of multivariate Cauchy distribution can be obtained using rmvc function of the LaplacesDemon package [32] in R.

Based on above methods, operators can fit the WPFE for all wind farms with a multivariate CD. Suppose that the vector $\tilde{\boldsymbol{p}}_t^w = (\tilde{p}_{1,t}^w, \tilde{p}_{2,t}^w, ..., \tilde{p}_{k,t}^w)^T$ represents the output of all $k$ wind farms in time period $t$. The fitting location vector is denoted as $\mu_t = (\mu_{1,t}, \mu_{2,t}, ..., \mu_{k,t})^T$ and the fitting scale matrix is denoted as $\Sigma_t$.

## III. MATHEMATICAL MODEL FORMULATION OF A-CCRTD

This section provides the formulation of A-CCRTD model. The declaration of the variables can be referred to the nomenclature section. Since the load forecast results are accurate enough with the state-of-the-art prediction technology [33], only the wind power prediction errors are considered in this paper. Actually, the model is scalable to incorporate the uncertainties of the load demand.

### A. Objective Function

$$F = \min \sum_{t=1}^{T} \left\{ \begin{array}{l} \sum_{i=1}^{N} CF_{i,t}(P_{i,t}^s) + \sum_{j=1}^{J} CF_{j,t}(P_{j,t}^a) + \\ \sum_{j=1}^{J} E[CR_{j,t}^+(\tilde{w}_t)] + \sum_{j=1}^{J} E[CR_{j,t}^-(\tilde{w}_t)] \end{array} \right\} \quad (8)$$

where $CF_{i,t}(\cdot)$ and $CF_{j,t}(\cdot)$ are the generation cost of AGC and non-AGC units, respectively. $CR_{j,t}^+(\cdot)$ and $CR_{j,t}^-(\cdot)$ represent the upward and downward regulation cost (corrective control cost) of AGC units, respectively; we can also regard these terms as the penalty cost of overestimation and underestimation of wind power output. The detailed formulations are listed below.

### 1) Generation cost

The generation costs of AGC and non-AGC units are expressed as quadratic functions of the power output:

$$CF_{i,t}(P_{i,t}^s) = a_{i,t}(P_{i,t}^s)^2 + b_{i,t}P_{i,t}^s + c_{i,t}$$
$$CF_{j,t}(P_{j,t}^a) = a_{j,t}(P_{j,t}^a)^2 + b_{j,t}P_{j,t}^a + c_{j,t} \quad (9)$$

### 2) Corrective control cost

Corrective control costs are caused by the mismatch between the actual wind power output $\tilde{w}_t$ and the arranged output $w_t$. The power mismatch should be balanced by the AGC units at any moment through certain principles. Considering the control principles in reality, the participation factors are usually assigned to AGC units proportional to their capacities. So the affine control strategy is established:

$$\tilde{p}_{j,t}^a = p_{j,t}^a - \alpha_j \cdot (\tilde{w}_t - w_t), \sum_{j=1}^{J} \alpha_j = 1(\alpha_j \geq 0) \quad (10)$$

where $\alpha_j$ is the participation factor of AGC unit $j$

The expectation of corrective costs is proportional to the expected positive and negative capacity deployed by the AGC units i.e.,

$$\begin{cases} E[CR_{j,t}^+(\tilde{w}_t)] = \gamma_{j,t}^+ \alpha_j \int_0^{w_t}(w_t - \tilde{\theta}_t)\varphi_t(\tilde{\theta}_t)d\tilde{\theta}_t \\ E[CR_{j,t}^-(\tilde{w}_t)] = \gamma_{j,t}^- \alpha_j \int_{w_t}^{\bar{w}_t}(\tilde{\theta}_t - w_t)\varphi_t(\tilde{\theta}_t)d\tilde{\theta}_t \end{cases} \quad (11)$$

where $\varphi_t(\cdot)$ is the *PDF* of the random variable in period $t$.

Assume that $\varphi_t(\tilde{w}_t)$ is the *PDF* of the summation of all wind farms' output in time period $t$. Referred to the mathematical properties of CD and the WPFE model in section II, we finally obtain:

$$\sum_{j=1}^{J} E[CR_{j,t}^+(\tilde{w}_t)] + \sum_{j=1}^{J} E[CR_{j,t}^-(\tilde{w}_t)] = \sum_{j=1}^{J} \left[ \begin{array}{l} A + B \cdot w_t - \dfrac{C \cdot \sqrt{\Sigma_{\tilde{w}_t}}}{2} \cdot \ln\left(1 + \left(\dfrac{w_t - \mu_{\tilde{w}_t}}{\sqrt{\Sigma_{\tilde{w}_t}}}\right)^2\right) \\ + C \cdot (w_t - \mu_{\tilde{w}_t})\arctan\dfrac{w_t - \mu_{\tilde{w}_t}}{\sqrt{\Sigma_{\tilde{w}_t}}} \end{array} \right] \quad (12)$$

where A, B and C are constants whose expressions are represented in Appendix, and $\mu_{\tilde{w}_t} = a_{\tilde{w}_t}^T \mu_t$, $\Sigma_{\tilde{w}_t} = a_{\tilde{w}_t}^T \Sigma_t a_{\tilde{w}_t}$. $a_{\tilde{w}_t}$ is a k-dimensional vector whose elements are all 1.

The objective function is convex and a relevant discussion is provided in Appendix.

### B. System constraints

Deterministic and chance constraints are shown below. For all constraints,

$$i \in \{1, 2, ..., n\}, j \in \{1, 2, ..., J\}, k \in \{1, 2, ..., K\},$$
$$l \in \{1, 2, ..., L\}, t \in \{1, 2, ..., T\}$$

$$\sum_{i=1}^{N} p_{i,t}^s + \sum_{j=1}^{J} p_{j,t}^a + \sum_{k=1}^{K} p_{k,t}^w = \sum_{d=1}^{D} p_{d,t}^d \quad (13)$$

$$\underline{P}_{j,t}^a \leq p_{j,t}^a \leq \bar{P}_{j,t}^a, \quad \underline{P}_{i,t}^s \leq p_{i,t}^s \leq \bar{P}_{i,t}^s, 0 \leq p_{k,t}^w \leq \bar{p}_{k,t}^w \quad (14)$$

$$\begin{cases} \Pr\{\alpha_j \cdot (w_t - \tilde{w}_t) + p_{j,t}^a \leq \bar{P}_{j,t}^a\} \geq 1 - \delta \\ \Pr\{\underline{P}_{j,t}^a \leq p_{j,t}^a + \alpha_j \cdot (\tilde{w}_t - w_t)\} \geq 1 - \delta \end{cases} \quad (15)$$

$$-RD_{i,t}^s \cdot \Delta T \leq p_{i,t}^s - p_{i,t-1}^s \leq RU_{i,t}^s \cdot \Delta T \quad (16)$$

$$\begin{cases} \Pr\{-RD_{j,t}^a \cdot \Delta T \leq \tilde{p}_{j,t}^a - \tilde{p}_{j,t-1}^a\} \geq 1 - \beta \\ \Pr\{\tilde{p}_{j,t}^a - \tilde{p}_{j,t-1}^a \leq RU_{j,t}^a \cdot \Delta T\} \geq 1 - \beta \end{cases} \quad (17)$$

$$\Pr\left\{R_t^+ \leq \sum_{j=1}^{J}(\bar{P}_{j,t}^a - \tilde{p}_{j,t}^a)\right\} \geq 1 - \varepsilon$$
$$\Pr\left\{R_t^- \leq \sum_{j=1}^{J}(\tilde{p}_{j,t}^a - \underline{P}_{j,t}^a)\right\} \geq 1 - \varepsilon \quad (18)$$

$$\Pr\left\{\left|\sum_{i=1}^{N} G_{l,i}p_{i,t}^s + \sum_{j=1}^{J} G_{l,j}\tilde{p}_{j,t}^a + \sum_{k=1}^{K} G_{l,k}\tilde{p}_{k,t}^w + \sum_{d=1}^{D} G_{l,d}p_{d,t}^d\right| \leq L_{l,t}\right\} \geq 1 - \eta \quad (19)$$

Where equation (13) is power balance constraint. Equation (14) is power generation limit constraint that means the scheduled power generation of AGC units, non-AGC units and wind farms cannot exceed their limits. Equation (15) is a chance constraint indicates that the actual regulation capacity of AGC units is guaranteed under a predefined confidence level. $\delta$ is the allowed probabilities of violation.

Equation (16) is the ramp-rate constraint for non-AGC units. Equation (17) limits the actual incremental output of AGC in adjacent time periods which is in the form of chance constraint. Obviously, unbalanced power caused by wind farm competes for ramp capability in real time operation, thus the additional power ramping requirement (APRR) should be included in dispatch process. $\beta$ is the pre-specified allowed probability of violation.

Equations (18) are reserve constraints to ensure system security under some contingency scenarios.

Equation (19) is transmission capacity constraint which represents the probability of transmission line overloading is no more than $\eta$. Note that the real-time unbalanced power allocated to each AGC units contributes to the active power on transmission lines, which is ignored in conventional CCED models.

### C. Solution Procedure

#### 1) Compact form and solution for chance constraints

Supposed that $A^{(g)}, B^{(g)} \in \mathbb{R}^n$, $D^{(g)} \in \mathbb{R}$, vector $\mathbf{u} \in \mathbb{R}^n$ represents decision variables and random vector $\tilde{\mathbf{y}} \in \mathbb{R}^n$, then chance constraints in A-CCRTD can be expressed in compact forms with equations (21) or (22):

$$\Pr\left[\left(A^{(g)}\right)^T \mathbf{u} + \left(B^{(g)}\right)^T \tilde{\mathbf{y}} \leq D^{(g)}\right] \geq 1 - \zeta^{(g)} \quad (20)$$

$$\Pr\left[\left(A^{(g)}\right)^T \mathbf{u} + \left(B^{(g)}\right)^T \tilde{\mathbf{y}} \geq D^{(g)}\right] \geq 1 - \zeta^{(g)} \quad (21)$$

According to the mathematical properties of CD described in section II, equations (20) and (21) are converted into deterministic constraints (22) and (23):

$$D^{(g)} - \left(A^{(g)}\right)^T \mathbf{u} \geq \left(B^{(g)}\right)^T \mu_{\tilde{y}} + \sqrt{\left(B^{(g)}\right)^T \Sigma_{\tilde{y}} B^{(g)}} \tan\left[\pi\left(1 - \zeta^{(g)} - \frac{1}{2}\right)\right] \quad (22)$$

$$D^{(g)} - \left(A^{(g)}\right)^T \mathbf{u} \leq \left(B^{(g)}\right)^T \mu_{\tilde{y}} + \sqrt{\left(B^{(g)}\right)^T \Sigma_{\tilde{y}} B^{(g)}} \tan\left[\pi\left(\zeta^{(g)} - \frac{1}{2}\right)\right] \quad (23)$$

All the chance constraints in A-CCRTD can be converted in the same way.

#### 2) The transformation of chance constraints in A-CCRTD

With WPFE model given in section II, chance constraint (15) is converted into constraint (24),

$$\begin{cases} \alpha_j w_t + p_{j,t}^a - \bar{P}_{j,t}^a \leq \alpha_j \cdot CDF_{\tilde{w}_t}^{-1}(\delta) \\ \alpha_j \cdot CDF_{\tilde{w}_t}^{-1}(1-\delta) \leq \alpha_j w_t + p_{j,t}^a - \underline{P}_{j,t}^a \end{cases} \quad (24)$$

where $CDF_{\tilde{w}_t}^{-1}(F) = \mu_{\tilde{w}_t} + \sqrt{\Sigma_{\tilde{w}_t}} \tan\left[\pi\left(F - \frac{1}{2}\right)\right]$.

Chance constraint (17) is converted into constraint (25):

$$\begin{cases} p_{j,t}^a - p_{j,t-1}^a + \alpha_j (w_t - w_{t-1}) - RU_{j,t}^a \cdot \Delta T \leq \alpha_j \cdot CDF_{w_{t,t-1}}^{-1}(\beta) \\ \alpha_j \cdot CDF_{w_{t,t-1}}^{-1}(1-\beta) \leq p_{j,t}^a - p_{j,t-1}^a + \alpha_j (w_t - w_{t-1}) + RD_{j,t}^a \cdot \Delta T \end{cases} \quad (25)$$

Where

$$CDF_{w_{t,t-1}}^{-1}(F) = a_{\tilde{w}_t} \mu_t - a_{\tilde{w}_{t-1}} \mu_{t-1} + \sqrt{a_{\tilde{w}_t}^T \Sigma_t a_{\tilde{w}_t} + a_{\tilde{w}_{t-1}}^T \Sigma_{t-1} a_{\tilde{w}_{t-1}}} \tan\left[\pi\left(F - \frac{1}{2}\right)\right]$$

Chance constraint (18) is converted into constraint (26):

$$\begin{cases} w_t + R_t^+ + \sum_{j=1}^{J} \left(p_{j,t}^a - \bar{P}_{j,t}^a\right) \leq CDF_{\tilde{w}_t}^{-1}(\varepsilon) \\ CDF_{\tilde{w}_t}^{-1}(1-\varepsilon) \leq w_t - R_t^- + \sum_{j=1}^{J} \left(p_{j,t}^a - \underline{P}_{j,t}^a\right) \end{cases} \quad (26)$$

And chance constraint (19) is converted into constraint (27):

$$\begin{cases} CDF_{a_l \tilde{P}_t^w}^{-1}(1-\eta) \leq L_l - \left[\sum_{i=1}^{N} G_{l,i} p_{i,t}^s + \sum_{j=1}^{J} G_{l,j} p_{j,t}^a + \left(\sum_{j=1}^{J} G_{l,j} \alpha_j\right) \sum_{k=1}^{K} p_{k,t}^w + \sum_{d=1}^{D} G_{l,d} p_{d,t}^d\right] \\ -L_l - \left[\sum_{i=1}^{N} G_{l,i} p_{i,t}^s + \sum_{j=1}^{J} G_{l,j} p_{j,t}^a + \left(\sum_{j=1}^{J} G_{l,j} \alpha_j\right) \sum_{k=1}^{K} p_{k,t}^w + \sum_{d=1}^{D} G_{l,d} p_{d,t}^d\right] \leq CDF_{a_l \tilde{P}_t^w}^{-1}(\eta) \end{cases} \quad (27)$$

where $CDF_{a_l \tilde{P}_t^w}^{-1}(F) = a_l^T \mu_t + \sqrt{a_l^T \Sigma_t a_l} \tan\left[\pi\left(F - \frac{1}{2}\right)\right]$, $\alpha_l$ is a K-dimensional vector whose $k$th element is $G_{l,k} - \left(\sum_{j=1}^{J} G_{l,j} \alpha_j\right)$.

Finally, the A-CCRTD model is analytically transformed to a convex objective function (represented by equations (9), (10), (12)) with deterministic linear constraints (including equations (13), (14), (16), (18), (24)-(27)).

**Remark**: The final model is convex and there is no approximation or iteration involved in transformation process. The fast computation performance of this model will be demonstrated in the numerical tests.

## IV. NUMERICAL TESTS

In this section, numerical tests were conducted to verify the effectiveness of the proposed method. First, the accuracy of CD in WPFE fitting was illustrated using real data of 20 wind farms in Southwest China. Then, the merits of the proposed model were justified on the modified IEEE 24-bus system. Meanwhile, the effects caused by the dependence between multiple wind farms in RTD was discussed. Finally, we verified the efficiency of A-CCRTD using the modified IEEE 118-bus system. Each real-time dispatch schedule was composed of 12 periods ($\Delta T$ = 5 min, T = 12), and only the solution of the first period was executed. The proposed model was solved using the IPOPT solver [35]. All simulations were implemented using Matlab R2015a on a laptop with an Intel Core i7 1.99 GHz processor and 8 GB of RAM.

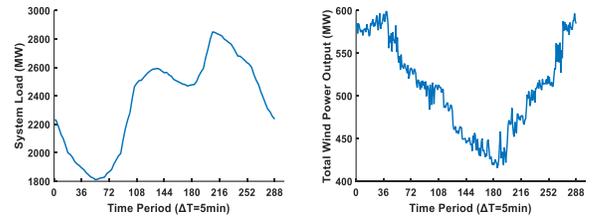

Fig.1. System loads and predicted total wind power output curve over 288 time periods.

Parameters of the modified IEEE 24-bus system are described below. The load profile of system is shown in the left of Fig. 1, where the valley-load periods are 01:30~07:30 and the peak-load periods are 15:30~21:30. The predicted profile of the total wind power output is shown in the right of Fig. 1, and the output feature obeys the rule that wind resources are greater during the nighttime. Four wind farms are connected at buses 7, 14, 16, and 21, which are denoted by #1, #2, #3, and #4, respectively. The capacity of the four wind farms is set to 240, 300, 80, and 180 MW, respectively. AGC units are connected at buses 5-8, 23, and 31-33, and the participation factor of each AGC unit is proportional to its capacity. If not specially specified, the price of upward AGC regulation capacity is \$12/MWh, and that of downward AGC regulation capacity is \$24/MWh. The normalized ramp rates of all units are defined with the ratio of the ramping capacities in MW/$\Delta T$ to their maximum capacities. The normalized ramp rates are assumed to be 0.05 and 0.1 for non-AGC units and AGC units, respectively. In addition, all confidence levels in this simulation are set to be 0.98. Details of the configuration and parameters of the modified IEEE 24-bus system and the modified IEEE 118-bus system can be accessed in [36].

*A. Comparison Fitting Acuracy of WPFE Using Different Distributions*

To illustrate the high accuracy of Cauchy distribution in WPFE fitting, 20 wind farms in the Southwest China with more than 80000 data were analyzed statistically. All actual and ultra-short-term forecasting data used here were retrieved from the electric power control center. We normalized all the forecast wind power and corresponding actual wind power within [0,1] [27] [28], and then fitted the WPFE with conditional distributions in different actual values. Most of the actual values belong to interval 0.0 p.u. to 0.7 p.u.. So we randomly selected two data sets here each containing around 7000 data pairs. In Data 1, the predicted wind power is concentrated near 0.1 p.u. and the Data 2 is concentrated around 0.4 p.u.. From Fig.2 and Fig.3, we concluded that the CD remarkably outperforms the other distributions, especially in capturing the kurtosis and tail behavior. Gaussian, Beta and Weibull distributions obviously underestimates the probability in the middle and overestimates the probability in the head/tail region. As far as the RMES listed in Table I, the CD is also much better than other distributions inaccuracy.

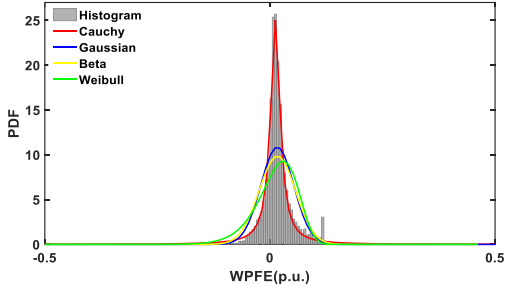

**Fig.2 PDF fitting results of different distributions using Data 1**

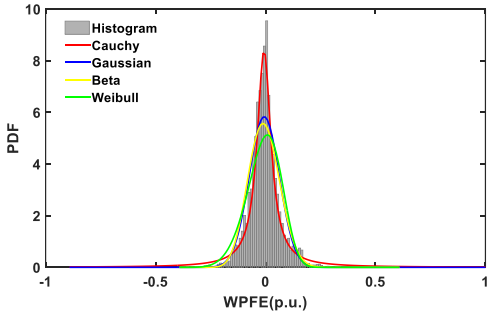

**Fig.3 PDF fitting results of different distributions using Data 2**

TABLE I
RMSEs OF DIFFERENT DISTRIBUTIONS

| Data Set | RMSE(p.u.) | | | |
|---|---|---|---|---|
| | Cauchy | Gaussian | Beta | Weibull |
| Data 1 | 0.3221 | 2.1144 | 2.2739 | 2.5273 |
| Data 2 | 0.3220 | 0.6365 | 0.7021 | 0.8695 |

*B. Comparison of A-CCRTD with conventional CCED*

In this subsection, CCED models without considering APRR and affine control strategy for AGC units were compared with the proposed A-CCRTD model. All simulations were for 21:00-22:00 and the cost was the sum of 12 dispatching periods. Monte Carlo simulations (MCS) with 10000 scenarios were conducted to compare the economic and security performance of A-CCRTD and others.

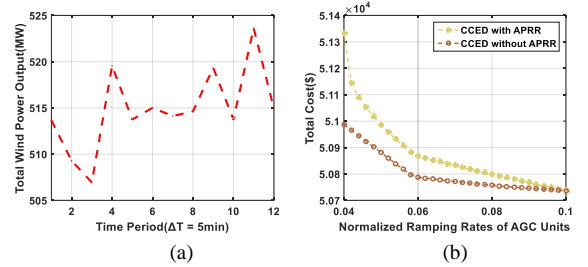

**Fig.4. (a) Predicted total wind power output curve in 21:00-22:00. (b) Total cost in two cases: CCED with APRR and CCED without APRR.**

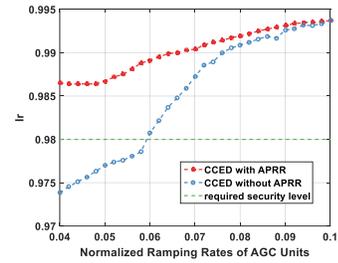

**Fig.5. Security level of unit ramping from time period 3 to time period 4.**

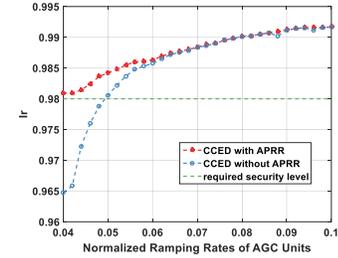

**Fig.6. Security level of unit ramping from time period 10 to time period 11.**

*1) The effect of APRR on ramping constraints*

In this simulation, the normalized ramp rates of all AGC units varied uniformly from 0.04 to 0.1. For the ease of comparison, security index for ramping resources is defined as

$$Ir = \frac{N_r}{N_M} \quad (28)$$

Where, $N_r$ is the average number of scenarios with sufficient ramping resources and $N_M$ is the total number of scenarios in MCS. Accordingly, a larger value of $Ir$ indicates a higher security level. Fig.4.(a) is the predicted profile of the total wind power output from 21:00 to 22:00. The total costs in each simulation are presented in Fig.3.(b), which shows that when the ramping rates of AGC units are low, taking into account the impact of APRR would increase the scheduling cost. The uneconomical scheduling results are compelling choice to prevent the lack of ramping resources in extreme scenarios. Fig.5 and Fig.6 are the Monte Carlo simulation results for AGC units in two different cases: (1) from time period 3 to time period 4; (2) from time period 10 to time period 11. Because of the rapid fluctuation of wind power in these two cases, which can be seen from Fig.4.(a), the security level of unit ramping without APRR cannot reach the required level in Fig.5 and Fig.6 due to the exhausted ramping resources when the ramping rates of AGC units are low. In addition, from Fig.5 and Fig.6, we can also conclude that only when there are abundant ramping resources in the system, the impact of APRR in ED can be ignored.

*2) The effect of affine control strategy on transmission capacity constraints*

In this simulation, we adjusted the transmission capacity of line #11 from 155MW to 170MW to illustrate the effect of affine control strategy on transmission capacity constraints. Similarly, we define security index for transmission capacity as

$$It = \frac{N_t}{N_M} \quad (29)$$

Where, $N_t$ is the average number of scenarios without transmission congestions for line #11 during all 12 dispatching periods. From Fig.7, it can be concluded that although the adoption of affine control strategy increases the operational cost, it guarantees a sufficient transmission capacity for security. This is because the redistribution of real-time power mismatch may cause transmission line congestion if we don't consider the control strategy of AGC in advance. i.e., it is necessary to incorporate regulation strategy of AGC units in schedule stage to prevent network congestions.

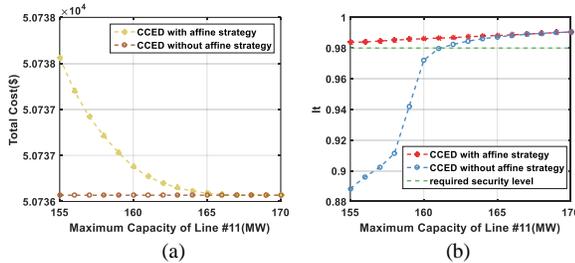

**Fig.7. Total cost(a) and security level(b) with different capacity of line #11. The cost is the summation of 12 dispatching periods and the security level is the average security level in 12 periods.**

### C. The effect of dependence between multi-wind farms

This experiment was carried out to test the influence of dependence between multiple wind farms on system performance. We used two different cases with and without considering the dependence between four wind farms for comparison: *Case* I: the location vector and scale matrix are consistent with the parameters mentioned before; *Case* II: the *PDF* of all wind farms are determined only by the marginal distribution respectively in *Case* I, that is to say the outputs of each wind farm are independent random variables. A-CCRTD with 12 periods was run in each cases. Based on the obtained schedule, 10000 random wind power scenarios were produced using the distribution parameters in *Case* I to illustrate the effect of dependence in terms of cost and risk level. From Table II, we can conclude that the consideration of dependence in ED reduces potential risk level at the expense of higher costs. This is because the integration of multi-wind farms amplifies the WPFE in our simulation. Therefore, the dependence of multi-wind farms should be considered in real-time power dispatch.

TABLE II
EFFECT OF CORRELATION ON ECONOMY AND OPERATIONAL RISK

| Case | Case I (dependent) | Case II (independent) |
|---|---|---|
| cost($) | 50736.09 | 50387.62 |
| Maximum risk level of reserve constraints | 1.57% | 2.38%(violate the predefined level) |
| Maximum risk level of unit ramping constraints | 0.83% | 0.78% |
| Maximum risk level on transmission line constraints | 0.73% | 0.64% |

### D. The Efficiency of A-CCRTD

The model size and computation times of A-CCRTD for both the IEEE 24-bus system and the IEEE 118-bus system are listed in table III. It is noteworthy that since the inverse *CDF* of CD is analytical, we can directly obtain the quantiles of chance constraints. In spite of 1225 variables and 6275 constraints involving in the A-CCRTD for the IEEE 118-bus system, it still can be solved in 7.23s. Therefore, this solution is sufficiently efficient for real-time application in large-scale power systems with high wind power penetration.

TABLE III
COMPUTATIONAL EFFICIENCY OF ACC-RTD MODEL

| System | 24 bus | 118 bus |
|---|---|---|
| CPU time (s) | 2.1611 | 7.2253 |
| Bus No. | 24 | 118 |
| Line No. | 38 | 181 |
| Unit and Wind farm No. | 33 | 79 |
| Variables No. | 625 | 1225 |
| Constraints No. | 2387 | 6275 |

## V. CONCLUSION

This paper proposes a A-CCRTD approach coordinating wind farms, non-AGC units, and AGC units. Two important features distinguish our model from the conventional CCED model. On the one hand, based on the accurate description of WPFE by Cauchy distribution and its promising mathematical properties, the A-CCRTD model is transformed equivalently to a convex optimization problem which is solved efficiently without any approximation. On the other hand, by incorporating affine control strategies of AGC units, our model takes into account both the APRR and the power variation on transmission lines caused by the allocation of real-time power mismatch between AGC units. Numerical tests demonstrate the superiority and rationality of the proposed approach compared with existing CCED models. Moreover, the importance of the dependence between multi-wind farms in real-time dispatch is also explored with Monte Carlo simulations. With practically acceptable computation effort even by general optimization solvers, the proposed approach makes the chance-constrained RTD practically tractable in real-time applications of large-scale power systems with high wind power penetration.

## VI. APPENDIX

### A. Constants in section III

$$B = -\frac{\alpha_j}{\pi}\left(\gamma_{j,t}^+ \arctan\frac{-\mu_{\tilde{w}_t}}{\sqrt{\Sigma_{\tilde{w}_t}}} + \gamma_{j,t}^- \arctan\frac{\bar{w}_t - \mu_{\tilde{w}_t}}{\sqrt{\Sigma_{\tilde{w}_t}}}\right) \quad C = \frac{\gamma_{j,t}^+ \alpha_j}{\pi} + \frac{\gamma_{j,t}^- \alpha_j}{\pi}$$

$$A = \frac{\alpha_j \cdot \sqrt{\Sigma_{\tilde{w}_t}}}{2\pi} \left( \gamma_{j,t}^+ \ln\left(1 + \left(\frac{\mu_{\tilde{w}_t}}{\sqrt{\Sigma_{\tilde{w}_t}}}\right)^2\right) + \gamma_{j,t}^- \ln\left(1 + \left(\frac{\bar{w}_t - \mu_{\tilde{w}_t}}{\sqrt{\Sigma_{\tilde{w}_t}}}\right)^2\right) \right)$$
$$+ \frac{\alpha_j \cdot \sqrt{\Sigma_{\tilde{w}_t}}}{\pi} \left( \gamma_{j,t}^+ \arctan \frac{-\mu_{\tilde{w}_t}}{\sqrt{\Sigma_{\tilde{w}_t}}} + \gamma_{j,t}^- \arctan \frac{\bar{w}_t - \mu_{\tilde{w}_t}}{\sqrt{\Sigma_{\tilde{w}_t}}} \right)$$

*B. The convexity of the objective function*

Suppose that $k_{1,j}$ and $k_{2,j}$ are the underestimation and overestimation cost coefficients, $v$ is the actual power output and $s$ is the scheduled decision variable. Then, the second order derivative of $\sum_{j=1}^{J} E\left[CR_{j,t}^+(v)\right] + \sum_{j=1}^{J} E\left[CR_{j,t}^-(v)\right]$ can be expressed as equation (30), where $p(v)$ is the *PDF* of $v$. Because $\alpha_j$, $k_{1,j}$, $k_{2,j}$ and $p(v)$ are all invariably positive, equation (12) is convex. Thus, the objective function is convex.

$$\frac{d^2}{ds^2}\left[\sum_{j=1}^{J}\left(k_{1,j} \cdot \alpha_j \int_s^{s_m}(v-s)p(v)dv + k_{2,j} \cdot \alpha_j \int_0^s(s-v)p(v)dv\right)\right]$$
$$= \frac{d}{ds}\left[\sum_{j=1}^{J}\left(k_{2,j}\alpha_j \int_0^s p(v)dv - k_{1,j}\alpha_j \int_s^{s_m}p(v)dv\right)\right] = \sum_{j=1}^{J}\left[\alpha_j(k_{1,j} + k_{2,j})p(v)\right] \quad (30)$$